# A note on some binomial sums

Johann Cigler
Fakultät für Mathematik
Universität Wien
A-1090 Wien, Nordbergstraße 15
Austria

johann.cigler@univie.ac.at

**Abstract**
Modifying an idea of E. Brietzke we give simple proofs for the recurrence relations of some sequences of binomial sums which have previously been obtained by other more complicated methods.

AMS Subject Classification: 05A10, 11B39.

*Key words*: Binomial coefficient, Fibonacci polynomial, Lucas polynomial, recurrence

**1. Introduction**
Modifying an idea in a recent paper by E. Brietzke [2] we give simple proofs for the recurrence relations of some sequences of binomial sums of the form

$$a(n) = \sum_{j \in \mathbb{Z}} z^j \binom{n}{\left\lfloor \frac{n-mj+k}{2} \right\rfloor},$$ which have previously been obtained by other methods.

In order to introduce the method we start with a well-known special case. Consider the sequence

$$a(n,5,k) = \sum_{j \in \mathbb{Z}} (-1)^j \binom{n}{\left\lfloor \frac{n-5j+k}{2} \right\rfloor} = (-1)^k \sum_{j \in \mathbb{Z}} t(n, k-5j) \text{ with } t(n,k) = (-1)^k \binom{n}{\left\lfloor \frac{n+k}{2} \right\rfloor}.$$

The recurrence of the binomial coefficients implies that $t(n,k) = -t(n-1, k-1) - t(n-1, k+1)$. For $n=0$ we get $t(0,0)=1, t(0,1)=-1$ and for all other $k \in \mathbb{Z}$ we have $t(0,k)=0$.

Let $s(n,k)$ on $\mathbb{N} \times \mathbb{Z}$ be the function which satisfies the same recurrence, but with initial values $s(0,k) = [k=0]$. Since $t(0,k) = s(0,k) - s(0,k-1)$ we get
$t(n,k) = s(n,k) - s(n,k-1)$.
Now consider the table

| $k$ | −3 | −2 | −1 | 0 | 1 | 2 | 3 |
|---|---|---|---|---|---|---|---|
| $s(0,k)$ | 0 | 0 | 0 | 1 | 0 | 0 | 0 |
| $s(1,k)$ | 0 | 0 | −1 | 0 | −1 | 0 | 0 |
| $s(2,k)$ | 0 | 1 | 0 | 2 | 0 | 1 | 0 |

Consider the vector space $\mathcal{F}$ of all functions $f(n,k)$ on $\mathbb{N}\times\mathbb{Z}$, such that for each $n\in\mathbb{N}$ there are only finitely many $k$ with $f(n,k)\neq 0$.

Let $T$ be the linear operator on $\mathcal{F}$ defined by $Ts(n,k) = s(n+2,k) - s(n+1,k) - s(n,k)$.
Then it is immediately verified that $Ts(0,k) = s(2,k) - s(1,k) - s(0,k) = 1$ for $|k|\leq 2$ and $Ts(0,k) = 0$ for $|k| > 2$.

This implies that $T\sum_{j\in\mathbb{Z}} s(0,k-5j) = 1$ for all $k\in\mathbb{Z}$ since $|k-5j|\leq 2$ for precisely one $j$.

Therefore $T\sum_{j\in\mathbb{Z}} t(0,k-5j) = T\sum_{j\in\mathbb{Z}} s(0,k-5j) - T\sum_{j\in\mathbb{Z}} s(0,k-1-5j) = 0$.

The same argument applies to $t(n,k)$, because each $t(n,k)$ is a linear combination of functions of the form $s(0,k+\ell) - s(0,k+\ell-1)$. This implies

$$a(n+2,5,k) - a(n+1,5,k) - a(n,5,k) = (-1)^k T\sum_{j\in\mathbb{Z}} (-1)^{5j+k}\left(\begin{array}{c} n \\ \left\lfloor\frac{n-5j+k}{2}\right\rfloor \end{array}\right) = 0.$$

Thus $a(n,5,k)$ satisfies the recurrence $a(n+2,5,k) - a(n+1,5,k) - a(n,5,k) = 0$.
Since the Fibonacci numbers $F_n$ satisfy the same recurrence with initial values $F_0 = 0$ and $F_1 = 1$, we get the following results:

For $k\equiv 0,1 \pmod{10}$ the initial conditions are $a(0,5,k) = a(1,5,k) = 1$ and therefore $a(n,5,k) = F_{n+1}$.

For $k\equiv 2,9 \pmod{10}$ we have $a(0,5,k) = 0, a(1,5,k) = 1$ and therefore $a(n,5,k) = F_n$.

For $k\equiv 3,8 \pmod{10}$ we get $a(0,5,k) = a(1,5,k) = 0$ and therefore $a(n,5,k) = 0$. Furthermore $a(n,5,k+5) = -a(n,5,k)$.

This result has first been proved by I. Schur [5] in a strengthened version. He has in fact given a $q$-analogue of this result. It has also been proved by G.E. Andrews [1].

## 2. A useful method

After this example let us consider a more general case.
For $a,b\in\mathbb{C}$ let $s(n,k,a,b)$ be the function on $\mathbb{N}\times\mathbb{Z}$ defined by $s(0,k,a,b) = [k=0]$ and the recurrence relation

$$s(n,k,a,b) = as(n-1,k-1,a,b) + bs(n-1,k,a,b) + as(n-1,k+1,a,b). \tag{1}$$

Let $\mathcal{F}_{a,b}$ be the subspace of $\mathcal{F}$ consisting of all functions satisfying the recurrence (1).

For any polynomial $p(x) = \sum_{i=0}^m a_i x^i$ we denote by $p(E)$ the linear operator on functions $f(n)$
defined by $p(E)f(n) = \sum_{i=0}^m a_i f(n+i)$. For $f(n,k)\in\mathcal{F}_{a,b}$ this means

$$p(E)f(n,k) = \sum_{i=0}^m a_i f(n+i,k).$$

We are looking for an operator $p(E)$ which has analogous properties as the operator $T$ in the above example.

To this end we define a sequence of polynomials $p_n(x,a,b) = \sum_{k=0}^{n} p_{n,k}(a,b)x^k$ by the recurrence

$$p_n(x,a,b) = (x-b)p_{n-1}(x,a,b) - a^2 p_{n-2}(x,a,b) \tag{2}$$

with initial values $p_0(x,a,b) = 1$ and $p_1(x,a,b) = x+a-b$.

**Lemma 1**

*For all $(m,k) \in \mathbb{N} \times \mathbb{Z}$ the following identity holds*

$$p_m(E,a,b)s(0,k,a,b) = \sum_{i=0}^{m} p_{m,i}(a,b)s(i,k,a,b) = a^m \left[|k| \leq m\right]. \tag{3}$$

**Proof**

It is immediately verified that (3) is true for $m=0$ and $m=1$.
If (3) has already been shown for $m-1$ and $m-2$ we get

$$p_m(E,a,b)s(0,k,a,b) = \left((E-b)p_{m-1}(E,a,b) - a^2 p_{m-2}(E,a,b)\right)s(0,k,a,b)$$

$$= p_{m-1}(E,a,b)s(1,k,a,b) - bp_{m-1}(E,a,b)s(0,k,a,b) - a^2 p_{m-2}(E,a,b)s(0,k,a,b)$$

$$= ap_{m-1}(E,a,b)s(0,k-1,a,b) + ap_{m-1}(E,a,b)s(0,k+1,a,b) - a^2 p_{m-2}(E,a,b)s(0,k,a,b)$$

$$= a^m \left[|k-1| \leq m-1\right] + a^m \left[|k+1| \leq m-1\right] - a^m \left[|k| \leq m-2\right] = a^m \left[|k| \leq m\right].$$

**Corollary**

*For each $k \in \mathbb{Z}$*

$$p_m(E,a,b) \sum_{j \in \mathbb{Z}} s(0, k-(2m+1)j, a, b) = \sum_{j \in \mathbb{Z}} p_m(E,a,b)s(0, k-(2m+1)j, a, b) = a^m. \tag{4}$$

**Proof**

For each $k \in \mathbb{Z}$ there is precisely one $j$ such that $|k-(2m+1)j| \leq m$.

As an easy application we consider for each $m \in \mathbb{N}$ the sequence

$$a(n, 2m+1, k) = \sum_{j \in \mathbb{Z}} (-1)^j \binom{n}{\left\lfloor \frac{n-(2m+1)j+k}{2} \right\rfloor}.$$

By the above argument it satisfies the recurrence relation $p_m(E,-1,0)a(n, 2m+1, k) = 0$.

Recall that the Fibonacci polynomials $F_n(x,s) = \sum_{k=0}^{n-1} \binom{n-1-k}{k} s^k x^{n-2k-1}$ are characterized by the recurrence $F_n(x,s) = xF_{n-1}(x,s) + sF_{n-2}(x,s)$ with initial conditions $F_0(x,s) = 0$ and $F_1(x,s) = 1$.
Now $p_n(x,-1,0)$ satisfies the recurrence $p_n(x,-1,0) = xp_{n-1}(x,-1,0) - p_{n-2}(x,-1,0)$ with initial values $p_0(x,-1,0) = 1$ and $p_1(x,-1,0) = x-1$.
Therefore $p_n(x,-1,0) = F_{n+1}(x,-1) - F_n(x,-1)$.

This gives

**Theorem 1**

*The sequence* $a(n, 2m+1, k) = \sum_{j \in \mathbb{Z}} (-1)^j \left( \left\lfloor \dfrac{n - (2m+1)j + k}{2} \right\rfloor \right)$ *satisfies the recurrence relation*

$$\left(F_{m+1}(E, -1) - F_m(E, -1)\right) a(n, 2m+1, k) = 0 \tag{5}$$

*for each* $k \in \mathbb{Z}$.

**Remark**

This theorem has been proved in [1] and [3] with more complicated methods.
It should be noted that $a(n, 2m+1, 0)$ has a simple combinatorial interpretation. It is the number of the set of all lattice paths in $\mathbb{R}^2$ which start at the origin, consist of $\left\lfloor \dfrac{n}{2} \right\rfloor$ northeast steps $(1,1)$ and $\left\lfloor \dfrac{n+1}{2} \right\rfloor$ southeast steps $(1,-1)$ and which are contained in the strip $-m-1 < y < m$. (cf. e.g. [4]).

## 3. A modification of the above method

The analogous result for the sequences $a(n, 2m, k)$ is a bit more complicated. Here we define a sequence of polynomials $q_n(x, a, b) = \sum_{k=0}^{n} q_{n,k}(a, b) x^k$ by the same recurrence

$$q_n(x, a, b) = (x - b) q_{n-1}(x, a, b) - a^2 q_{n-2}(x, a, b), \tag{6}$$

but with initial values $q_0(x, a, b) = 2$ and $q_1(x, a, b) = x - b$.

**Lemma 2**

*For all* $(m, k) \in \mathbb{N} \times \mathbb{Z}$ *the following identity holds*

$$q_m(E, a, b) s(0, k, a, b) = \sum_{i=0}^{m} q_{m,i}(a, b) s(i, k, a, b) = a^m \left[ |k| = m \right]. \tag{7}$$

**Proof**

(7) is true for $m = 1$ and $m = 2$ by inspection.
If it is already shown for $m - 1$ and $m - 2$ we get
$q_m(E, a, b) s(0, k, a, b) = \left((E - b) q_{m-1}(E, a, b) - a^2 q_{m-2}(E, a, b)\right) s(0, k, a, b)$
$= a q_{m-1}(E, a, b) s(0, k-1, a, b) + q_{m-1}(E, a, b) s(0, k+1, a, b) - a^2 q_{m-2}(E, a, b) s(0, k, a, b)$
$= a^m \left[ |k-1| = m-1 \right] + a^m \left[ |k+1| = m-1 \right] - a^m \left[ |k| = m-2 \right] = a^m \left[ |k| = m \right].$

As an application let $v(n,k) = \left(\left\lfloor \dfrac{n}{\frac{n+k}{2}} \right\rfloor\right)$.

Then $v(n,k) = v(n-1,k-1) + v(n-1,k+1)$ and $v(0,k) = [k \in \{0,1\}]$. Therefore
$v(n,k) = s(n,k,1,0) + s(n,k-1,1,0)$.

We have

$$a(n,2m,k) = \sum_{j \in \mathbb{Z}} (-1)^j \left(\left\lfloor \dfrac{n}{\frac{n-(2m)j+k}{2}} \right\rfloor\right) = \sum_{j \in \mathbb{Z}} \left(\left(\left\lfloor \dfrac{n}{\frac{n-(2m)2j+k}{2}} \right\rfloor\right) - \left(\left\lfloor \dfrac{n}{\frac{n-(2m)2j+k-2m}{2}} \right\rfloor\right)\right)$$

$$= \sum_{j \in \mathbb{Z}} \left(s(n,k-4mj,1,0) - s(n,k-2m-4mj,1,0)\right) + \sum_{j \in \mathbb{Z}} \left(s(n,k-1-4mj,1,0) - s(n,k-1-2m-4mj,1,0)\right).$$

Here we get $q_m(E,1,0) \sum_{j \in \mathbb{Z}} \left(s(n,k-4mj,1,0) - s(n,k-2m-4mj,1,0)\right) = 0$,

because $s(n,k-4mj,1,0) - s(n,k-2m-4mj,1,0)$ is a linear combination of functions of the form $s(0,k+\ell-4mj,1,0) - s(0,k+\ell-2m-4mj,1,0)$, where for $k+\ell-4mj = m$ we get $k+\ell-4mj-2m = -m$ and the sums cancel. For other values $k$ the sum vanishes.

Recall that the Lucas polynomials $L_n(x,s) = \sum_{k=0}^{n-1} \binom{n-k}{k} \dfrac{n}{n-k} s^k x^{n-2k}$ are characterized by the recurrence $L_n(x,s) = xL_{n-1}(x,s) + sL_{n-2}(x,s)$ with initial conditions $L_0(x,s) = 2$ and $L_1(x,s) = x$. Therefore $p_n(x,-1,0) = L_n(x,-1)$ and we obtain

**Theorem 2**

*The sequence* $a(n,2m,k) = \sum_{j \in \mathbb{Z}} (-1)^j \left(\left\lfloor \dfrac{n}{\frac{n-(2m)j+k}{2}} \right\rfloor\right)$ *satisfies the recurrence relation*

$$L_m(E,-1)a(n,2m,k) = 0. \tag{8}$$

The same method can be applied to the sum

$$a(n,m,k,z) = \sum_{j \in \mathbb{Z}} z^j \left(\left\lfloor \dfrac{n}{\frac{n-mj+k}{2}} \right\rfloor\right) = \sum_{j \in \mathbb{Z}} z^{2j} \left(\left\lfloor \dfrac{n}{\frac{n-2mj+k}{2}} \right\rfloor\right) + \sum_{j \in \mathbb{Z}} z^{2j-1} \left(\left\lfloor \dfrac{n}{\frac{n-2mj+k+m}{2}} \right\rfloor\right).$$

Here we get
$$L_m(E,-1)a(0,m,k,z) = L_m(E,-1)\sum_{j \in \mathbb{Z}} z^{2j} v(0,k-2mj) + L_m(E,-1) \sum_{j \in \mathbb{Z}} z^{2j-1} v(0,k+m-2mj).$$

In this case we have $L_m(E,-1)v(0,k-2mj) = \begin{cases} 1 & \text{if } k = 2mj - m \\ 1 & \text{if } k = 2mj + m \\ 0 & \text{else} \end{cases}$

or $L_m(E,-1)v(0,k-2mj) = v(0,k-m-2mj) + v(0,k+m-2mj)$. This implies

$$L_m(E,-1)a(0,m,k,z) = \sum_{j\in\mathbb{Z}} z^{2j}(v(0,k-m-2mj)+v(0,k+m-2mj))$$

$$+\sum_{j\in\mathbb{Z}} z^{2j-1}(v(0,k+2m-2mj)+v(0,k-2mj)) = \left(z+\frac{1}{z}\right)a(0,m,k,z).$$

Thus we get

$$\left(L_m(E,-1)-\left(z+\frac{1}{z}\right)\right)a(0,m,k,z) = 0.$$

This implies in turn

**Theorem 3**

*The sequence* $a(n,m,k,z) = \sum_{j\in\mathbb{Z}} z^j \binom{n}{\left\lfloor\frac{n-mj+k}{2}\right\rfloor}$ *satisfies the recurrence relation*

$$\left(L_m(E,-1)-\left(z+\frac{1}{z}\right)\right)a(n,m,k,z) = 0. \tag{9}$$

Theorems 2 and 3 have previously been proved in [3] with other methods.